\documentclass[11pt,twoside]{article}
\usepackage[hmargin=0.8in,vmargin=0.9in]{geometry}
\usepackage{fancyhdr}
\usepackage{graphicx}
\usepackage{subfigure}
\usepackage{amssymb}
\usepackage{amsmath}
\usepackage{amsfonts}
\usepackage{theorem}
\usepackage{mathrsfs}
\usepackage{mathtools}
\usepackage{bm}
\usepackage{color}
\usepackage{setspace}
\usepackage{exscale}
\usepackage{relsize}
\usepackage{epstopdf}
\usepackage{float}
\usepackage{cite}

\usepackage{color}

\usepackage{booktabs,multirow} % for much better looking tables
\usepackage{array} % for better arrays (eeg matrices) in maths
\usepackage{paralist} % very flexible & customisable lists (eeg. enumerate/itemize, etc.)
\usepackage{verbatim} % adds environment for commenting out blocks of text & for better verbatim
\usepackage{subfigure} % make it possible to include more than one captioned figure/table in a single float
\theoremstyle{plain}			% use "default" font
\newtheorem{thm}{Theorem}[section]

\newtheorem{prop}[thm]{Proposition}

{\theorembodyfont{\rmfamily}}

\newenvironment{DA}{{\flushleft \bf Declarations:}}{}
\newtheorem{rmk}{Remark}[section]

\allowdisplaybreaks[1]

\numberwithin{equation}{section}
\numberwithin{figure}{section}
\numberwithin{table}{section}

\newcommand\eref[1]{(\ref{#1})}

\newcommand*\xbar[1]{%
  \hbox{%
    \vbox{%
      \hrule height 0.5pt % The actual bar
      \kern0.4ex%         % Distance between bar and symbol
      \hbox{%
        \kern-0.05em%      % Shortening on the left side
        \ensuremath{#1}%
        \kern-0.00em%      % Shortening on the right side
      }%
    }%
  }%
}

\newcommand{\bmF}{\bm{\mathcal{F}}}

\newcommand{\bmu}{\bm{u}}

\newcommand{\mf}{\bm{f}}

\newcommand{\mo}{\bm{0}}

\newcommand{\dt}{\Delta t}
\newcommand{\dx}{\Delta x}

\newcommand{\dxi}{\Delta \xi}

\newcommand{\hf}{{\frac{1}{2}}}

\newcommand{\jph}{{j+\frac{1}{2}}}
\newcommand{\jmh}{{j-\frac{1}{2}}}

\def\softd{{\leavevmode\setbox1=\hbox{d}%
		\hbox to 1.05\wd1{d\kern-0.4ex{\char039}\hss}}}

\title{Numerical Analysis of Stabilization for Random Hyperbolic Systems of Conservation Laws}
\author{Shaoshuai Chu\thanks{Department of Mathematics, RWTH Aachen University, Aachen, 52056, Germany; {\tt chu@igpm.rwth-aachen.de}},
~Michael Herty\thanks{Department of Mathematics, RWTH Aachen University, Aachen, 52056, Germany; Department of Mathematics and Applied
Mathematics, University of Pretoria, Private Bag X20, Hatfield 0028, South Africa; {\tt herty@igpm.rwth-aachen.de}}, and Alexander Kurganov
\thanks{Department of Mathematics and Shenzhen International Center for Mathematics, Southern University of Science and Technology,
Shenzhen, 518055, China; {\tt alexander@sustech.edu.cn}}}

\begin{document}
\date{}
\maketitle
\begin{abstract}
This paper extends the deterministic Lyapunov-based stabilization framework to random hyperbolic systems of conservation laws,
where uncertainties arise in boundary controls and initial data. Building on the finite-volume discretization method from [{\sc M. Banda and
M. Herty}, Math. Control Relat. Fields, 3 (2013), pp. 121--142], we introduce a stochastic discrete Lyapunov function to prove the
exponential decay of numerical solutions for systems with random perturbations. For linear systems, we derive explicit decay rates, which
depend on boundary control parameters, grid resolutions, and the statistical properties of the random inputs. Theoretical decay rates are
verified through numerical examples, including boundary stabilization of the linear wave equations and linearized shallow-water flows with
random perturbations. We also present the decay rates for a nonlinear example and for the linearized Saint-Venant system with source terms.
\end{abstract}



\section{Introduction}
The study of stabilization for hyperbolic systems has been an active field of research, particularly in the context of control theory and
numerical analysis. A key aspect of the stabilization analysis in hyperbolic partial differential equations (PDEs) involves feedback
boundary conditions, which have been extensively studied in various contexts; see, e.g., \cite{coron2007control,Bastin2016}. These studies
highlight the importance of well-posed boundary conditions in ensuring controllability and stabilizability. The interplay between analytical
techniques and numerical discretization methods has also been explored in \cite{Bardos1979}, where finite-volume schemes are employed to
preserve the stability properties of the continuous problem. Building on these insights, later studies have broadened the stabilization
toolbox by developing approaches based on Lyapunov functions, feedback boundary control, and further refined numerical discretization
techniques; see, e.g.,
\cite{bastin2007lyapunov,coron2007control,coron2008dissipative,coron2007strict,li2009exact,li2008controllabilite,li2003exact,li2005exact,xu2002exponential,banda2013numerical}. In particular, for the stabilization of Saint-Venant system of shallow-water equations,
substantial progress has been achieved in the analysis of related models that incorporate source terms. For instance, dissipative boundary
conditions for one-dimensional (1-D) nonlinear hyperbolic systems were investigated in \cite{coron2008dissipative}, while a Lyapunov-based
framework for exponential stability of 1-D linear hyperbolic balance laws was developed in \cite{DBC2012}. Extensions to more complex
systems, such as the Saint-Venant-Exner model, were studied in \cite{DDTK2017}, where backstepping-based stabilization techniques were
proposed. More recently, numerical boundary feedback stabilization for non-uniform hyperbolic balance laws was introduced in \cite{BW2020}.
These works highlight the importance of accounting for source terms and provide valuable insights into the stabilization of the Saint-Venant
system and related models.

Recent studies have introduced Lyapunov-based approaches to establish exponential stability for deterministic hyperbolic balance laws. These approaches provide a rigorous framework for analyzing the long-term behavior of solutions in the presence of well-posed boundary conditions; see, e.g., \cite{Krstic2008}. The application of discrete Lyapunov functions has been particularly successful in proving the stability of numerical schemes, as it allows for explicit stability bounds and decay rates under controlled settings; see, e.g., \cite{banda2013numerical, Bastin2016,Meurer2013}.

This paper aims to bridge this gap by extending the Lyapunov-based stabilization framework introduced in \cite{banda2013numerical} 
to hyperbolic balance laws with stochastic perturbations.  This extension is motivated by the increasing interest in stochastic PDEs and their applications in modeling uncertainties in physical systems, including fluid dynamics, traffic flow, and networked systems, where random variations must be accounted for in the stability analysis; see, e.g., \cite{Karatzas1991, Fleming1975, Holden2015,Ghanem2003}. While previous studies have extensively addressed deterministic stabilization techniques, the extension to random systems remains relatively underdeveloped. The incorporation of randomness into the stabilization framework introduces new mathematical challenges, particularly in the construction of appropriate Lyapunov functions that account for stochastic perturbations. Similar approaches were developed in the context of stochastic control theory, where Lyapunov techniques were adapted to random dynamical systems; see, e.g., \cite{Kushner1967}. By incorporating randomness into the framework established in \cite{banda2013numerical}, this study advances the theoretical understanding of stabilization for stochastic hyperbolic balance laws and provides computational results that support the theoretical findings.

The rest of the paper is organized as follows. \S\ref{sec2} provides a brief overview of the studied hyperbolic systems, outlining the
governing equations and boundary conditions. In \S\ref{sec3},  we introduce the numerical discretization scheme and present theoretical
stability results using the proposed stochastic Lyapunov framework. \S\ref{sec4} discusses generalizations of the stabilization framework,
including modifications of boundary conditions and extensions to broader classes of systems. In \S\ref{sec5}, we validate the theoretical
findings through numerical experiments, including the stabilization of linearized shallow-water equations under random perturbations.
We also demonstrate the decay rates in a nonlinear shallow-water example. Moreover, we extend the analysis to random hyperbolic
systems of balance laws with boundary damping, for which both theoretical results and numerical validations are presented; see Appendix C.
Finally, we give some concluding remarks and potential directions for future research in \S\ref{sec6}.

\section{An Overview of the Studied Hyperbolic Systems}\label{sec2}
In this section, we give a brief introduction to the studied problem and refer the reader to \cite{banda2013numerical} for more details.

We consider the following initial-boundary value problems for 1-D nonlinear hyperbolic PDEs:
\begin{equation}\label{1.1}
\bmu_t+ \mf (\bmu)_x=\mo, \quad x\in[0,1],\,\,t\in(0,T],
\end{equation}
with the initial conditions
\begin{equation}\label{1.1a}
\bmu(0,x,\xi)=\bm \psi(x,\xi),\quad x\in[0,1],
\end{equation}
and the boundary conditions, which will be specified later. In \eref{1.1} and \eref{1.1a},  $t$ is the time, $x$ is the spatial variable, $\xi$  is the random variable, $\bmu=\big(u^{(1)},\ldots, u^{(p)}\big)^\top$ is the vector of unknowns, and $\mf:\mathbb{R}^p \to \mathbb{R}^p$ denotes a possibly nonlinear smooth flux function. The system \eref{1.1} is assumed to be strictly hyperbolic, meaning that the Jacobian matrix $F(u) ={\partial \mf}/{\partial \bmu}$ has real distinct eigenvalues. For smooth solutions $\bmu(x,t,\xi)$,  \eref{1.1} can be rewritten in a quasi-linear form as
\begin{equation}\label{1.2}
\bmu_t + F(\bmu) \bmu_x = 0, \quad \bmu(0, x, \xi) = \bm \psi(x,\xi).
\end{equation}

The eigenvalues of $F(\bmu)$ are denoted by $\Lambda_i$, $i = 1, \dots, p$, and we assume that there exists $m$ such that for any $\bmu$,
$$
\Lambda_p <  \ldots< \Lambda_{m+1} <0 < \Lambda_{m} < \ldots < \Lambda_{1}.
$$
This classification allows us to define the decomposition of $\bm u$ into components associated with positive and negative eigenvalues:
$$
\bmu =
\begin{pmatrix}
    \bmu^+ \\
    \bmu^-
\end{pmatrix},
$$
where $\bmu^+ \in \mathbb{R}^m$ and $\bmu^- \in \mathbb{R}^{p-m}$ correspond to the positively and negatively propagating characteristics, respectively. Similarly, we decompose the flux function:
$$
F(\bmu)=\begin{pmatrix}F^+(\bmu)\\F^-(\bmu)\end{pmatrix}, 
$$
where $F^+:\mathbb R^p\to\mathbb R^{m\times p}$, $F^-:\mathbb R^p\to\mathbb R^{(p-m)\times p}$. In addition, $G:\mathbb R^p\to\mathbb R^p$
is a possibly nonlinear boundary operator
$$
G(\bmu)=\begin{pmatrix}G^+(\bmu)\\G^-(\bmu)\end{pmatrix},
$$
which determines how information propagates at the domain boundaries with $G^+:\mathbb R^p\to\mathbb R^m$, and
$G^-:\mathbb R^p\to\mathbb R^{p-m}$ representing the corresponding flux components.

As in \cite{banda2013numerical}, we adopt a general feedback-type boundary condition, which reads as
$$
\begin{pmatrix}\bmu^+(t,0,\xi)\\\bmu^-(t,1,\xi)\end{pmatrix}=G\begin{pmatrix}\bmu^+(t,1,\xi)\\\bmu^-(t,0,\xi)\end{pmatrix}. 
$$

\section{The Discretization Scheme} \label{sec3}
In this section, we study the numerical schemes for boundary $L^2$-stabilization of 1-D nonlinear random hyperbolic systems. We first consider a simple case with $m \equiv p$, that is, $\Lambda_i >0$, $\forall i$,  $F(\bmu)$ is diagonal, and $G$ is a linear operator with $G_{ij}=0$, $i \neq j$.

\begin{prop}\label{prop1}
Assume that $F$ is diagonal for all $\bmu\in B_\varepsilon(0)  \subset \mathbb{R}^p$ where $B_\varepsilon(0)$ is an open ball in $\mathbb{R}^p$ centered at the origin and with radius $\varepsilon$. Also assume that 
\begin{equation}\label{1.5}
F(\bmu) ={\rm diag}\big( \Lambda_1(\bmu), \ldots, \Lambda_m(\bmu)\big), \quad \Lambda_i>0, \quad \Lambda_i(\bmu) \ne \Lambda_j(\bmu), \quad \forall i \ne j,
\end{equation}
and the boundary conditions are prescribed as 
\begin{equation}\label{1.6}
\bmu(t,0,\xi)=\mathcal{K} \bmu(t,1,\xi),
\end{equation}
where 
\begin{equation*}
 \mathcal{K}={\rm diag}\big(\kappa_1, \ldots, \kappa_m\big), \quad \kappa_i>0, \,\, \forall i.
 \end{equation*}
 If $\max\limits_{i=1,\ldots,m} \kappa_i <1$, then the equilibrium $\bmu \equiv 0$ for \eref{1.2} with \eref{1.5} and \eref{1.6} is exponentially stable. 
\end{prop}
The proof of Proposition \ref{prop1}  immediately follows from the fact that  $\rho(\mathcal{K})	=\max\limits_{i=1,\ldots,m} \kappa_i<1$.

We now introduce the numerical discretization. We assume that $\delta>0$ is sufficiently small such that  $M_\delta(0) \subset B_\varepsilon(0)$ with
$$ 
M_\delta(0) := \{ \bmu : |u^{(i)}| \leq \delta, \,\, i = 1, \dots, m \}.
$$

We introduce a uniform grid $(x_j,\,\xi_k)$ for $j=0,\ldots,M$ and $k=0,\ldots,K$ with $x_{j+1}-x_{j}\equiv\dx$, $\xi_{k+1}-\xi_{k}\equiv\dxi$. The temporal grid is chosen such that the CFL condition holds
\begin{equation}\label{1.6a}
 \lambda \frac{\dt}{\dx} \le 1, \quad \lambda:=\max\limits_{i=1,\ldots,m} \max\limits_{\bmu\in M_{\delta(0)}} \Lambda_i (\bmu)
\end{equation}
and $t^n=n \dt$ for $n=0,1,\ldots N$, where by possibly further reducing $\dt$, one can enforce $N \dt=T$.

Assume that the initial conditions $\bmu^{\,0}_{j,k}=\bm \psi(x_j,\xi_k)$, $j=1,\ldots,M$, $k=0,\ldots,K$ are available. We then use the boundary conditions \eref{1.6} to obtain 
$u^{(i),0}_{0,k}=\kappa_i u^{(i),0}_{M,k}$, $k=0,\ldots,K$, and then evolve the solution in time by the upwind scheme used in \cite{banda2013numerical}:

\begin{align}
u_{j,k}^{(i),n+1} &= u_{j,k}^{(i),n} - \frac{\Delta t}{\Delta x} \Lambda_i \left( \bmu_{j,k}^{n} \right) \left( u_{j,k}^{(i),n} - u_{j-1,k}^{(i),n} \right),\quad j=1,\ldots,M,\,k=0,\ldots,K,\,i=1,\ldots,m, \label{1.10} \\
u_{0,k}^{(i),n+1} &= \kappa_i u_{M,k}^{(i),n+1}, \quad k=0,\ldots,K,\,i=1,\ldots,m.\label{1.10a}
\end{align}

Notice that the scheme \eref{1.10}--\eref{1.10a} is nonconservative, but this is not a problem since in this paper we only consider smooth solutions of \eref{1.1}--\eref{1.1a}.  It needs to be noted that if ${u}^{(i),n}_{j,k} \in M_{\delta(0)}$ then 
\begin{equation}\label{1.10b}
0 < D^{(i),n}_{j,k}:=\frac{\dt }{\dx} \Lambda_i\big(\bm{u}^{n}_{j,k}\big) \le \frac{\Lambda_i \big(\bm{u}^{n}_{j,k}\big)}{\lambda} \le 1,
\end{equation}
and the scheme \eref{1.10}--\eref{1.10a} can be rewritten as
\begin{equation}\label{1.10c}
u^{(i),n+1}_{j,k}=u^{(i),n}_{j,k}\big(1-D^{(i),n}_{j,k}\big)+u^{(i),n}_{j-1,k}D^{(i),n}_{j,k}, \quad u^{(i),n+1}_{0,k} = \kappa_i u^{(i),n+1}_{M,k}.
\end{equation}

\begin{prop}\label{prop2}
Assume \eref{1.5}, \eref{1.6}, $\big|u^{(i),0}_{j,k}\big| \le \delta$ for $i=1,\ldots,m$, $j=1, \ldots, M$,  $k=0,\ldots, K$, and $\max\limits_{i=1,\ldots,m} \kappa_i \le 1$. Then the scheme \eref{1.10}--\eref{1.10a} satisfies
\begin{equation*}
 \big|u^{(i),n}_{j,k}\big|  \le \delta, \quad \forall\, i,j,k,n.
\end{equation*}
\end{prop}

{\emph {Proof}}. We prove this proposition by induction. Assume that $\big|u^{(i),n}_{j,k} \big|\le \delta$ and therefore $D^{(i),n}_{j,k} \in(0, 1]$.  Then, \eref{1.10c} implies
\begin{equation*}
\begin{aligned}
\big|u^{(i),n+1}_{j,k}\big|& = \big|u^{(i),n}_{j,k} \big(1-D^{(i),n}_{j,k}\big)+u^{(i),n}_{j-1,k} D^{(i),n}_{j,k}\big| 
	               \le \big|u^{(i),n}_{j,k}\big|\big(1-D^{(i),n}_{j,k}\big)+\big|u^{(i),n}_{j-1,k}\big| D^{(i),n}_{j,k}
	                \le \delta,
\end{aligned}
\end{equation*}
for all $i=1,\ldots,m$, $j=1,\ldots,M$, and $k=0,\ldots,K$. In addition, \eref{1.10a} gives $\big|u^{(i),n+1}_{0,k}\big| = \kappa_i \big|u^{(i),n+1}_{M,k}\big| \le \delta$, which completes the proof of the proposition.

$\hfill\square$

Using Proposition \ref{prop2} and \eref{1.10b}, one immediately obtains that for all $i$, 
\begin{equation}\label{1.13}
 \max\limits_{j,k,n}D^{(i),n}_{j,k} \le \frac{\dt}{\dx} \max\limits_{\bm{u}\in M_\delta(0)}  \Lambda_i(\bm{u})=:D^{\rm max}_{i} \le 1;
\end{equation}
and 
\begin{equation}\label{1.14}
\frac{\dt}{\dx}\Lambda_i( \bm \eta) \ge \frac{\dt}{\dx} \min\limits_{\bmu\in M_\delta(0)} \Lambda_i(\bm{u})=:D^{\rm min}_i >0, \quad \forall \, {\bm \eta}\in M_\delta(0).
\end{equation}

Next, we obtain a bound on the discrete first spatial derivative, which is an expected result in view of the analytical bound on $\bmu_x$; see e.g., \cite{coron2007control}.
\begin{prop}\label{prop3}
Assume \eref{1.5}, \eref{1.6}, $\big|u^{(i),0}_{j,k}\big| \le \delta$, 
\begin{equation}\label{1.14a}
 \bigg|\dfrac{u^{(i),0}_{j,k}-u^{(i),0}_{j-1,k}}{\dx}\bigg| \le \delta <1, 
\end{equation}
for $i=1,\ldots,m$, $j=1,\ldots,M$,  $k=0, \ldots, K$, and $\max\limits_{i=1,\ldots,m} \kappa_i \le 1$. Further assume that 
\begin{equation*}
\bigg|\dfrac{u^{(i),n}_{1,k}-u^{(i),n}_{0,k}}{\dx}\bigg| \le \delta e^{n \dt J^{(i)}_{\rm max}}, \quad   
\end{equation*}
with 
\begin{equation*}
J^{(i)}_{\rm max}:=\max\limits_{\bm \eta \in M_\delta(0)}|| \nabla_{\bmu}\Lambda_i(\bm{\eta}) ||_\infty,
\end{equation*}
for $i=1,\ldots,m$, $k=0, \ldots, K$, and $n=0,\ldots,N$. Then the scheme \eref{1.10}--\eref{1.10a} satisfies
\begin{equation}\label{1.15}
\bigg|\dfrac{u^{(i),n}_{j,k}-u^{(i),n}_{j-1,k}}{\dx}\bigg|\le \delta e^{n \dt J^{(i)}_{\rm max}},
\end{equation}
for all $i=1,\ldots,m$, $k=0, \ldots, K$, $n=0,\ldots,N$, and $j \ge 2$.
\end{prop}

{\emph {Proof}}. We prove this proposition by induction. First, we write expressions for $u^{(i),1}_{j,k}$ and $u^{(i),1}_{j-1,k}$ using \eref{1.10}, subtract them, divide by $\dx$, and obtain that for $j\ge2$, the following holds:
$$
\begin{aligned}
\bigg|\frac{u^{(i),1}_{j,k}-u^{(i),1}_{j-1,k}}{\dx}\bigg| & \le (1-D^{(i),0}_{j,k}) 	\bigg|\frac{u^{(i),0}_{j,k}-u^{(i),0}_{j-1,k}}{\dx}\bigg|+D^{(i),0}_{j-1,k} 	\bigg|\frac{u^{(i),0}_{j-1,k}-u^{(i),0}_{j-2,k}}{\dx}\bigg|\\
		& \overset{\eref{1.14a}}{\le} \delta \Big(1-\Big(D^{(i),0}_{j,k}-D^{(i),0}_{j-1,k}\Big)\Big) 
		 \overset{\eref{1.10b}}{\le} \delta \Big(1+\frac{\dt}{\dx}\big|\Lambda_i\big(\bmu^{0}_{j,k}\big)-\Lambda_i\big(\bmu^{0}_{j-1,k}\big) \big|\Big)\\
        & \overset{\eref{1.14a}}{\le} \delta\big (1+ \dt J^{(i)}_{\rm max} \big)
		 \le \delta e^{\dt J^{(i)}_{\rm max}}. 
\end{aligned}
$$
Assume now that \eref{1.15} holds at all time levels until $t=t^n$. Repeating the previous computations for the difference between $u^{(i),n+1}_{j,k}$ and $u^{(i),n}_{j,k}$ for $j\ge 2$, results in 
\begin{equation*}
\begin{aligned}
  \bigg|\frac{u^{(i),n+1}_{j,k}-u^{(i),n+1}_{j-1,k}}{\dx}\bigg|& \overset{\eref{1.10}, \eref{1.14a}}{\le}\delta e^{n \dt J^{(i)}_{\rm max}} \Big(1-\Big(D^{(i),n}_{j,k}-D^{(i),n}_{j-1,k}\Big)\Big)\\
    \hspace{3cm} & \le  \ldots \le \delta e^{n \dt J^{(i)}_{\rm max}} e^{ \dt J^{(i)}_{\rm max}}
    =\delta e^{(n+1) \dt J^{(i)}_{\rm max}},
\end{aligned}
\end{equation*}
which completes the proof of the proposition.

$\hfill\square$

We now introduce the discrete Lyapunov function at time $t=t^n$ with positive coefficients $\mu_i$, $i=1,\ldots,m$ as 
\begin{equation}\label{1.16}
{\cal{L}}^n = \dx \dxi \sum_{i=1}^{m}  \sum_{j=1}^{M} \sum_{k=1}^{K} \big(u^{(i),n}_{j,k}\big)^2 e^{-\mu_ix_j} \rho(\xi_k),
\end{equation}	
where $\rho(\xi)$ is the probability density function. As in the continuous case, the Lyapunov function can be bounded by 
\begin{equation}\label{1.17}
\min\limits_{i=1,\ldots,m}e^{-\mu_i x_{M}} \dx \dxi\sum_{i=1}^{m} \sum_{j=1}^{M}  \sum_{k=1}^{K} \big(u^{(i),n}_{j,k}\big)^2 \rho(\xi_k) \le {\cal L}^n \le  \dx \dxi \sum_{i=1}^{m} \sum_{j=1}^{M}  \sum_{k=1}^{K} \big(u^{(i),n}_{j,k}\big)^2  \rho(\xi_k).
\end{equation}

\begin{thm}\label{thm1}
Assume that \eref{1.5} and \eref{1.6} hold,  then if $\kappa_i$ and $\mu_i$ are bounded by 
\begin{equation}\label{1.18} 
 0 < \kappa_i <\sqrt{\frac{D^{\rm min}_i}{D^{\rm max}_i}},\qquad  \mu_i  \le \ln \Bigg(\bigg(\sqrt{\frac{D^{\rm max}_i}{D^{\rm min}_i}}\kappa_i\bigg)^{-2}\Bigg),
\end{equation}
for $i=1,\ldots,m$, and the following bounds on the initial and boundary values:
\begin{equation*}
\big|u^{(i),0}_{j,k}\big| \le \delta, \quad \bigg|\dfrac{u^{(i),0}_{j,k}-u^{(i),0}_{j-1,k}}{\dx}\bigg| \le \delta, \quad \bigg|\dfrac{u^{(i),n}_{1,k}-u^{(i),n}_{0,k}}{\dx}\bigg|\le \delta e^{n \dt J^{(i)}_{\rm max}},
\end{equation*}
are satisfied  for $i=1,\ldots,m$, $j=1,\ldots, M$,  and $k=0,\ldots,K$ with 
\begin{equation}\label{1.18c}
 \delta  = \min\Bigg\{ 1, \varepsilon, \frac{\dx}{2\dt}\min\limits_{i=1,\ldots,m} \Bigg\{  \frac{\mu_i D^{\rm min}_i e^{-T J^{(i)}_{\rm max}}}{J^{(i)}_{\rm max}}   \Bigg\}  \Bigg\}. 
\end{equation}
Then, for the numerical solution $u^{(i),n}_{j,k}$ defined by \eref{1.10}--\eref{1.10a}, the Lyapunov function \eref{1.16} satisfies 
\begin{equation}\label{1.19}
{\cal L}^n \le e^{-\nu t^n} {\cal L}^0, \quad n=0,1,\ldots,N,
\end{equation}
with 
\begin{equation}\label{1.19a}
 \nu =\frac{\dx}{2\dt} \min\limits_{i=1,\ldots,m} \Big\{  \mu_iD^{\rm min}_i   e^{-\mu_i \dx}   \Big\}.
\end{equation}
Moreover, $u^{(i),n}_{j,k}$ is exponentially stable in the discrete $L^2$-norm, namely,  
\begin{equation}\label{1.20}
\dx \dxi \sum_{i=1}^{m} \sum_{j=1}^{M}  \sum_{k=1}^{K} \big(u^{(i),n}_{j,k}\big)^2 \rho(\xi_k)  \le \max\limits_{i=1,\ldots,m}e^{\mu_ix_M-\nu t^n} \dx \dxi \sum_{i=1}^{m} \sum_{j=1}^{M}  \sum_{k=1}^{K} \big(u^{(i),0}_{j,k}\big)^2 \rho(\xi_k),
\end{equation}
for $n=0,\ldots,N$. 

\end{thm}

{\emph {Proof}}. First, we use \eref{1.10c} and an inequality $2ab \le a^2 +b^2$ to show that 
\begin{equation*}
\big(u^{(i),n+1}_{j,k}\big)^2 \le  \big(u^{(i),n}_{j,k}\big)^2 \big(1-D^{(i),n}_{j,k}\big) +\big(u^{(i),n}_{j-1,k}\big)^2 D^{(i),n}_{j,k},
\end{equation*}
and hence
\begin{equation}\label{1.21}
\big(u^{(i),n+1}_{j,k}\big)^2 - \big(u^{(i),n}_{j,k}\big)^2  \le D^{(i),n}_{j,k}  \Big(  \big(u^{(i),n}_{j-1,k}\big)^2 - \big(u^{(i),n}_{j,k}\big)^2 \Big).
\end{equation}
We then use the definition \eref{1.16} and obtain 
\begin{align*}
\frac{{\cal L}^{n+1}-{\cal L}^n}{\dt }& \overset{\eref{1.21}}{\le} \frac{\dx\dxi}{\dt}\sum_{i=1}^{m}\sum_{j=1}^{M}\sum_{k=1}^{K} D^{(i),n}_{j,k}\Big(\big(u^{(i),n}_{j-1,k}\big)^2- \big(u^{(i),n}_{j,k}\big)^2\Big)e^{-\mu_i x_j}\rho(\xi_k) \\
& \le \frac{\dx\dxi}{\dt} \Bigg[ \sum_{i=1}^{m} \sum_{j=1}^{M-1}\sum_{k=1}^{K}D^{(i),n}_{j+1,k}\big(u^{(i),n}_{j,k}\big)^2e^{-\mu_i x_{j+1}}\rho(\xi_k) 
+\sum_{i=1}^{m}\sum_{k=1}^{K}D^{(i),n}_{1,k}\big(u^{(i),n}_{0,k}\big)^2 e^{-\mu_i x_1}\rho(\xi_k)\\
&\hspace{1.5cm}-\sum_{i=1}^{m}\sum_{j=1}^{M}\sum_{k=1}^{K} D^{(i),n}_{j,k} \big(u^{(i),n}_{j,k}\big)^2 e^{-\mu_i x_j}\rho(\xi_k)
\Bigg]\\
&\overset{\eref{1.10a}}{\le} \frac{\dx \dxi}{\dt}\Bigg[ \sum_{i=1}^{m}\sum_{j=1}^{M-1}\sum_{k=1}^{K}D^{(i),n}_{j+1,k}\big(u^{(i),n}_{j,k}\big)^2e^{-\mu_i x_{j}} e^{-\mu_i \dx}\rho(\xi_k)\\
& \hspace{1.7cm}+\sum_{i=1}^{m}\sum_{k=1}^{K}D^{(i),n}_{1,k}\kappa_i^2\big(u^{(i),n}_{M,k}\big)^2 e^{-\mu_i x_1}\rho(\xi_k)
-\sum_{i=1}^{m}\sum_{j=1}^{M}\sum_{k=1}^{K} D^{(i),n}_{j,k} \big(u^{(i),n}_{j,k}\big)^2 e^{-\mu_i x_j}\rho(\xi_k)\Bigg].
\end{align*}
Next, we define a ghost value 
\begin{equation}\label{1.22}
D^{(i),n}_{M+1,k}:=D^{(i),n}_{M,k}, 
\end{equation}
add and subtract $\dfrac{\dx\dxi}{\dt}\sum\limits_{i=1}^{m}\sum\limits_{k=1}^{K} D^{(i),n}_{M+1,k}\big(u^{(i),n}_{M,k}\big)^2e^{-\mu_i(x_{M}+\dx)}\rho(\xi_k)$ to the right-hand side (RHS) of the previous inequality, which results in 
\begin{equation}\label{1.23}
\begin{aligned}
\frac{{\cal L}^{n+1}-{\cal L}^n}{\dt }&  \le  \frac{\dx \dxi}{\dt}  \Bigg[ \sum_{i=1}^{m} \sum_{j=1}^{M} \sum_{k=1}^{K} \big(u^{(i),n}_{j,k}\big)^2 \big(D^{(i),n}_{j+1,k}e^{-\mu_i\dx}-D^{(i),n}_{j,k}\big)e^{-\mu_i x_j} \rho(\xi_k) \\
&+ \sum_{i=1}^{m} \sum_{k=1}^{K} \bigg\{ D^{(i),n}_{1,k}  \kappa^2_i  e^{-\mu_i x_1}
-  D^{(i),n}_{M+1,k} e^{-\mu_i(x_{M}+\dx)} \bigg\}\big(u^{(i),n}_{M,k}\big)^2 \rho(\xi_k)\Bigg].
\end{aligned}
\end{equation}
According to the mean value theorem, there exists $\bm \eta \in {\rm conv }(\bm{u}^{\,n}_{j,k}, \bm{u}^{\,n}_{j+1,k})$ satisfying 
\begin{equation}\label{1.24}
\begin{aligned}
D^{(i),n}_{j+1,k}&=D^{(i),n}_{j,k}+\frac{\dt}{\dx}\nabla_{\bm{u}}\Lambda_i(\bm{\eta})(\bm{u}^{\,n}_{j+1,k}-\bm{u}^{\,n}_{j,k}) =D^{(i),n}_{j,k}+\dt \sum_{\ell=1}^{m}\partial_{\bmu_\ell} \Lambda_i({\eta^{(\ell)}}) \frac{u^{(\ell),n}_{j+1,k}  -u^{(\ell),n}_{j,k} }{\dx}\\
& \overset{ \eref{1.15}}{\le} D^{(i),n}_{j,k}+\dt \sum_{\ell=1}^{m} J^{(i)}_{\rm max} \delta e^{T J^{(i)}_{\rm max}} \overset{ \eref{1.18c}}{\le} D^{(i),n}_{j,k}+\dfrac{\mu_i \dx}{2}D^{\rm min}_i, 
\end{aligned}
\end{equation}
for $i=1,\ldots,m$, $k=0,\ldots,K$, and $j=1,\ldots,M$ (for $j=M$, one needs to also use \eref{1.22} to establish the last inequality in \eref{1.24}).

Therefore, one obtains
\begin{equation}\label{1.25}
\begin{aligned}
D^{(i),n}_{j+1,k} e^{-\mu_i \dx}-D^{(i),n}_{j,k}& \overset{ \eref{1.24}}{\le} \Big(  D^{(i),n}_{j,k} + \frac{\mu_i \dx}{2}D^{\rm min}_i\Big) e^{-\mu_i \dx}-D^{(i),n}_{j,k}  \\
& \overset{ \eref{1.14}}{\le} D^{\rm min}_i \big(e^{-\mu_i \dx}-1\big) + \frac{\mu_i\dx}{2} D^{\rm min}_i e^{-\mu_i \dx}\le -\frac{\mu_i \dx}{2} D^{\rm min}_i e^{-\mu_i \dx}.
\end{aligned}
\end{equation}
Next, using \eref{1.25} in \eref{1.23} leads to 
\begin{equation}\label{1.26}
\begin{aligned}
\frac{{\cal L}^{n+1}-{\cal L}^n}{\dt }&  \le  -\dfrac{\dx \dxi}{2 \dt} \sum_{i=1}^{m} \sum_{j=1}^{M}   \sum_{k=1}^{K} D^{\rm min}_i e^{-\mu_i \dx} \mu_i \dx e^{-\mu_i x_j} \big(u^{(i),n}_{j,k}\big)^2 \rho(\xi_k)\\
& + \frac{\dx \dxi}{\dt} \sum_{i=1}^{m} \sum_{k=1}^{K} \Big(D^{(i),n}_{1,k} \kappa_i^2 e^{-\mu_i x_1} -D^{(i),n}_{M+1,k}e^{-\mu_i(x_{M}+\dx)}\Big )\big(u^{(i),n}_{M,k}\big)^2 \rho(\xi_k).
\end{aligned}
\end{equation}
The last term on the RHS of \eref{1.26} is nonpositive since 
\begin{equation*}
\begin{aligned}
     & D^{(i),n}_{1,k}\kappa_i^2e^{-\mu_i x_1} -D^{(i),n}_{M+1,k}D^{(i),n}_{M+1,k}e^{-\mu_i(x_{M}+\dx)}
 \overset{ \eref{1.13}}{\le}   D^{\rm max}_i \kappa_i^2 e^{-\mu_i x_1} -D^{\rm min}_i e^{-\mu_i(x_{M}+\dx)}  \\
 =   & D^{\rm min}_i e^{-\mu_i (x_{M}+\dx)} \Bigg(\Bigg(\sqrt{\frac{D^{\rm max}_i}{D^{\rm min}_i}}       \kappa_i   \Bigg)^2 e^{\mu_i (x_{M}+\dx-x_1)} -1              \Bigg)\\
 =   & D^{\rm min}_i e^{-\mu_i (x_{M}+\dx)}\Bigg(\Bigg(\sqrt{\frac{D^{\rm max}_i}{D^{\rm min}_i}}         \kappa_i \Bigg)^2 e^{\mu_i} -1              \Bigg) \overset{ \eref{1.18}}{\le} 0.
\end{aligned}
\end{equation*}
Hence, 
\begin{equation*}
\frac{{\cal L}^{n+1}-{\cal L}^n}{\dt } \le -\dfrac{\dx \dxi}{2\dt} \sum_{i=1}^{m} \sum_{j=1}^{M}  \sum_{k=1}^{K}D^{\rm min}_i e^{-\mu_i \dx} \mu_i \dx e^{-\mu_i x_j} \big(u^{(i),n}_{j,k}\big)^2  \rho(\xi_k)
\overset{ \eref{1.19a}}{\le} - \nu {\cal L}^n.
\end{equation*}
Therefore, we obtain 
\begin{equation*}
{\cal L}^{n+1} - {\cal L}^n \le -\nu \dt {\cal L}^n,
\end{equation*}
and by recursively applying the previous inequality, we end up with  
\begin{equation*}
{\cal L}^{n} \le (1-\nu \dt)^{n} {\cal L}^0 \le e^{-\nu t^n} {\cal L}^0 = e^{-\nu t^{n}}{\cal L}^0.
\end{equation*}
Finally, the inequality \eref{1.20} follows from \eref{1.19} using the definition \eref{1.16} and  estimate \eref{1.17}.

$\hfill\square$
\begin{rmk}
From the proof of Theorem \ref{thm1}, we observe that its statement is also true in the case
$-\sqrt{\frac{D^{\min}_i}{D^{\max}_i}}<\kappa_i<0$, and we omit the proof for the sake of brevity. The case $\kappa_i=0$ is a special case.
Here, the boundary condition for the component $i$ is equal to zero and this state is propagated through the domain. From the practical
point of view, this is not interesting since the system would be controlled by precisely the state (namely zero), which it is
intended to reach. The previous proof yields in the case $\kappa_i=0$ no bound on the component $\mu_i$. In the case $\kappa_i=0$ for all
$i=1,\ldots,m$, we obtain exponential convergence for any rate $\nu>0$.  
\end{rmk}
\begin{rmk}
Note that the results are independent of the choice of either a general or uniform probability density function since the decay rate $\nu$
is independent of the values of the probability density $\rho(\xi)$ function. 
\end{rmk}

\section{Extensions of Theorem \ref{thm1}}\label{sec4}
In this section, we discuss some extensions and modifications of Theorem \ref{thm1}.

\subsection{Extension of Theorem \ref{thm1} to the Case of Linear $F$}
In the linear case, a stronger result can be obtained. Assuming
\begin{equation}\label{1.28}
F(u)={\rm diag}(\Lambda_i)^m_{i=1}, \, \Lambda_i >0, \, \Lambda_i \ne \Lambda_j,\, i \ne j
\end{equation}
holds, then \eref{1.6a}, \eref{1.13}, and \eref{1.14}, yield
\begin{equation*}
\lambda= \max\limits_{i=1,\ldots,m} \Lambda_i, \,\,\, D^{(i),n}_{j,k} = D^{\rm max}_i=D^{\rm min}_i = \frac{\dt}{\dx}\Lambda_i.
\end{equation*}
Note that in this case, $D^{(i),n}_{j,k}$ is independent of $j$, $k$, and $n$. Therefore, in the proof of Theorem \ref{thm1}, the estimate \eref{1.15} on the discrete derivative of $u^{(i),n}_{j,k}$ is not needed. Hence, one can use only Proposition \ref{prop1} to derive the following estimate as in the proof of Theorem \ref{thm1}:
\begin{equation*}
\frac{{\cal L}^{n+1}-{\cal L}^n}{\dt } \le -\dfrac{\dx}{\dt}\min\limits_{i=1,\ldots,m} \big\{ \mu_i D^{\rm min}_i  e^{-\mu_i \dx}\big\} {\cal L}^n.
\end{equation*}
Observe that compared with \eref{1.19a} an improved bound on the decay rate is obtained. Note that \eref{1.18} reduces now to 
\begin{equation*}
0 < \kappa_i <1, \quad \mu_i \le \ln(\kappa_i^{-2}).
\end{equation*}

As one can see, all bounds are now independent of $\delta$, and hence a finite terminal time $T$ is not needed. The results of the discussion above are summarized in the following Theorem.
\begin{thm}\label{thm4.1}
Assume that \eref{1.28} holds. For any $\kappa_i$, $i=1,\ldots,m$ such that $0 < \kappa_i <1$, there exists $\mu_i>0$, $i=1,\ldots,m$ such that for all initial data $\bm \psi_{j,k}$ the numerical solution $\bm u^n_{j,k}$ defined by \eref{1.10}--\eref{1.10a} satisfies 
$$
{\cal L}^n \le e^{-\nu t^n} {\cal L}^0, \,\, n=0,\ldots,N,
$$
for $\nu =  \min\limits_{i=1,\ldots,m}\big\{  \Lambda_i \mu_i e^{-\mu_i \dx}   \big\}$ with $\mu_i  \le \ln(\kappa^{-2}_i)$.   Furthermore, $\bm u^n_{j,k}$ is exponentially stable in the discrete $L^2$-norm and \eref{1.20} holds.
\end{thm}

\subsection{Extension of Theorem \ref{thm1} to the Case $m \ne p$}\label{sec4.2}
Recall that in \S \ref{sec3}, we assumed that $m=p$, that is, all transport coefficients $\Lambda_i(\cdot)$ were positive. When $m\neq p$, we consider a different Lyapunov function for the negative transport coefficients. Let
\begin{equation}\label{1.31}
{\cal L}^n = \dx \dxi \sum_{i=1}^{m} \sum_{j=1}^{M} \sum_{k=1}^{K} e^{-\mu_i x_j} \big(u^{(i),n}_{j,k}\big)^2 \rho(\xi_k)+\dx \dxi \sum_{i=m+1}^{p} \sum_{j=1}^{M} \sum_{k=1}^{K} e^{\mu_i x_j} \big(u^{(i),n}_{j,k}\big)^2 \rho(\xi_k),
\end{equation}
where, as before, $\mu_i>0$, $\forall i$, and it is assumed that for all $\bm\eta \in M_{\delta(0)}$
\begin{equation}\label{1.32}
\Lambda_i (\bm \eta)>0, \,\, i=1,\ldots,m, \,\, {\rm and}\,\,\, \Lambda_i(\bm \eta)<0, \, i=m+1,\ldots,p.
\end{equation}
Furthermore, the discrete boundary conditions need to be modified accordingly, namely, we set 
\begin{equation}\label{1.32a}
u^{(i),n}_{0,k}= \kappa_i u^{(i),n}_{M,k}, \,\, i=1,\ldots,m, \,\, k=0,\ldots K,
\end{equation}
and 
\begin{equation}\label{1.33}
u^{(i),n}_{M+1,k}=\kappa_i u^{(i),n}_{1,k},\, \, i=m+1,\ldots,p, \,\, k=0,\ldots K.
\end{equation}
The upwind scheme is also different now. While \eref{1.10}--\eref{1.10a} still can be used for $i\le m$, for $i=m+1,\ldots,p$, the upwinding direction changes and thus \eref{1.10}--\eref{1.10a} should be replaced with 
\begin{align}
\hspace{-0.3cm} u^{(i),n+1}_{j,k}& =u^{(i),n}_{j,k}-\frac{\dt}{\dx}\Lambda_i(\bm{u}^{\,n}_{j+1,k}) \big(u^{(i),n}_{j+1,k}-u^{(i),n}_{j,k}\big), \,\, j=1,\ldots,M,\,k=0,\ldots,K,\,i=m+1,\ldots,p,\label{1.34}\\[1.ex]
\hspace{-0.3cm} u^{(i),n+1}_{M+1,k}& = \kappa_i u^{(i),n+1}_{1,k}, \quad  k=0,\ldots,K,\,i=m+1,\ldots,p.\label{1.34a}
\end{align}
Similarly, the CFL condition \eref{1.6a} and the definitions of $D^{(i),n}_{j,k}$, $D^{\rm max}_i$, and $D^{\rm min}_i$ in \eref{1.10b}, \eref{1.13}, and \eref{1.14} should be adjusted to 
\begin{equation}\label{1.38}
\lambda \frac{\dt}{\dx} \le 1, \,\, \lambda:=\max\limits_{i=1,\ldots,p} \max\limits_{\bm{u}\in M_{\delta(0)}}|\Lambda_i (\bm u)|,
\end{equation}
and
\begin{equation}\label{1.40}
  0 < D^{(i),n}_{j,k}:= \frac{\dt}{\dx} |\Lambda_i(\bm{u}^{\,n}_{j,k})| \le 1, \quad  D^{\rm max}_{i}:=\frac{\dt}{\dx} \max\limits_{\bm{u}\in M_\delta(0)}  |\Lambda_i(\bm{u})|, \quad  D^{\rm min}_i:=\frac{\dt}{\dx} \min\limits_{\bmu\in M_\delta(0)} |\Lambda_i(\bm{u})|,
\end{equation}
respectively.

We now consider $i\ge m+1$, and rewrite \eref{1.34} as 
\begin{equation*}
u^{(i),n+1}_{j,k} = u^{(i),n}_{j,k}+\frac{\dt}{\dx}|\Lambda_i(\bm{u}^{\,n}_{j+1,k})| \big(u^{(i),n}_{j+1,k}-u^{(i),n}_{j,k}\big) 
		         \overset{\eref{1.40}}{=}  u^{(i),n}_{j,k} \big(1-D^{(i),n}_{j+1,k}\big)+D^{(i),n}_{j+1,k}u^{(i),n}_{j+1,k},
\end{equation*}
This leads to the following estimate:
\begin{equation}\label{1.40a}
\big(u^{(i),n+1}_{j,k}\big)^2-\big(u^{(i),n}_{j,k}\big)^2  \le D^{(i),n}_{j+1,k} \Big( \big(u^{(i),n}_{j+1,k}\big)^2-\big(u^{(i),n}_{j,k}\big)^2   \Big),
\end{equation}
which can be derived precisely the same way \eref{1.21} was obtained. We then multiply \eref{1.40a} by $e^{\mu_i x_j}$ and sum over $j$ to obtain 
\begin{equation}\label{1.40b}
\begin{aligned}
    & \hspace{-0.8cm} \sum_{j=1}^{M}   \Big( \big(u^{(i),n+1}_{j,k}\big)^2-\big(u^{(i),n}_{j,k}\big)^2   \Big) e^{\mu_i x_j}
\le  \sum_{j=1}^{M}  D^{(i),n}_{j+1,k}\Big( \big(u^{(i),n}_{j+1,k}\big)^2-\big(u^{(i),n}_{j,k}\big)^2 \Big)e^{\mu_i x_j}\\
  & \hspace{-1.3cm}= \sum_{j=1}^{M}   \big(u^{(i),n}_{j,k} \big)^2 \big(D^{(i),n}_{j,k} e^{-\mu_i \dx}-D^{(i),n}_{j+1,k} \big) e^{\mu_i x_j} 
    -  D^{(i),n}_{1,k}\big(u^{(i),n}_{1,k}\big)^2e^{\mu_ix_0} +D^{(i),n}_{M,k}\big(u^{(i),n}_{M+1,k}\big)^2 e^{\mu_i x_{M}}\\
    \hspace{-0.2cm}\overset{\eref{1.22},\eref{1.33}}{=}&  \sum_{j=1}^{M}   \big(u^{(i),n}_{j,k}\big)^2 \big(D^{(i),n}_{j,k}e^{-\mu_i \dx}-D^{(i),n}_{j+1,k} \big)e^{\mu_i x_j}
   +\big(u^{(i),n}_{1,k}\big)^2\big( D^{(i),n}_{M,k}\kappa_i^2 e^{\mu_i x_{M}}-D^{(i),n}_{1,k}e^{\mu_ix_0}\Big).
\end{aligned}
\end{equation}
Next, assuming similar bounds on the discrete derivative of the initial and boundary values as in Proposition \ref{prop3}, the first term on the RHS of \eref{1.40b} can be bounded using 
\begin{equation*}
D^{(i),n}_{j,k} e^{-\mu_i \dx} -D^{(i),n}_{j+1,k} \le  \frac{\mu_i \dx}{2} D^{\rm min}_i. 
\end{equation*}
The second term on the RHS of \eref{1.40b} can be bounded as before using
$$
\begin{aligned}
    & D^{(i),n}_{M,k}\kappa_i^2e^{\mu_i x_{M}}-D^{(i),n}_{1,k}e^{\mu_ix_0} 
\le  D^{\rm max}_i \kappa_i^2 e^{\mu_i x_{M}}-D^{\rm min}_ie^{\mu_ix_0}
   = D^{\rm min}_i e^{\mu_ix_0}\bigg( \bigg( \sqrt{\frac{D^{\rm max}_i}{D^{\rm min}_i}}  \kappa_i    \bigg)^2 e^{\mu_i}-1              \bigg) \overset{\eref{1.18}}{\le} 0.
\end{aligned}
$$
Combining the last two estimates and proceeding as in the proof of Theorem \ref{thm1}, we obtain 
\begin{equation*}
\begin{aligned}
    &\frac{\dx\dxi}{\dt}\sum_{i=m+1}^{p}\sum_{j=1}^{M}\sum_{k=1}^{K} \Big(\big(u^{(i),n+1}_{j,k}\big)-\big(u^{(i),n}_{j,k}\big)^2  \Big)e^{\mu_i x_j} \rho(\xi_k)\\
\le &-\frac{\dx\dxi}{2\dt}\sum_{i=m+1}^{p}\sum_{j=1}^{M}\sum_{k=1}^{K} D^{\rm min}_i   e^{-\mu_i \dx} \mu_i \dx e^{\mu_i x_j}\big(u^{(i),n}_{j,k}\big)^2 \rho(\xi_k) \\
\le & -\nu \dx 	\dxi\sum_{i=m+1}^{p}\sum_{j=1}^{M}\sum_{k=1}^{K} e^{\mu_i x_j}\big(u^{(i),n}_{j,k}\big)^2 \rho(\xi_k),
\end{aligned}
\end{equation*}
where $\nu$ is defined by 
\begin{equation*}
 \nu =\frac{\dx}{2\dt} \min\limits_{i=m+1,\ldots,p} \Big\{  \mu_iD^{\rm min}_i   e^{-\mu_i \dx}   \Big\}.
\end{equation*}
Hence, under suitable assumptions on the discrete gradients, there is an exponential decay without restriction on the sign of $\Lambda_i$.

We summarize the obtained results  in the following theorem.
\begin{thm}\label{thm3}
Assume that \eref{1.18}--\eref{1.18c}, \eref{1.32},  and \eref{1.38} hold. We also assume that 
$$\bigg|\dfrac{u^{(i),n}_{M+1,k}-u^{(i),n}_{M,k}}{\dx}\bigg|\le \delta e^{n \dt J^{(i)}_{\rm max}}$$
is valid for $i=m+1,\ldots,p$.
Then, for the  numerical solution $u^{(i),n}_{j,k}$ defined by \eref{1.10}--\eref{1.10a} for $i=1,\ldots,m$ and by \eref{1.34}--\eref{1.34a} for $i=m+1,\ldots,p$, the Lyapunov function \eref{1.31} satisfies \eref{1.19} for $\nu$ given by   
\begin{equation}\label{1.40d}
 \nu =\frac{\dx}{2\dt} \min\limits_{i=1,\ldots,p} \Big\{  \mu_iD^{\rm min}_i   e^{-\mu_i \dx}   \Big\}.
\end{equation}
Moreover, $u^{(i),n}_{j,k}$ is exponentially stable in the discrete $L^2$-norm, that is, \eref{1.20} is valid.
\end{thm}

\subsection{Extension of Theorem \ref{thm1} to Different Boundary Conditions}
Let us consider a specific case with $m=1$, $p=2$, so that the Lyapunov function \eref{1.31} reduces to 
\begin{equation}\label{1.41}
{\mathcal L}^n = \dx \dxi \sum_{j=1}^{M} \sum_{k=1}^{K} \Big[ e^{-\mu_1 x_j} \big(u^{(1),n}_{j,k}\big)^2 \rho(\xi_k)+ e^{\mu_2 x_j} \big(u^{(2),n}_{j,k}\big)^2 \rho(\xi_k) \Big],
\end{equation}
and $u^{(1),n}_{j,k}$ and $u^{(2),n}_{j,k}$ are computed by \eref{1.10}--\eref{1.10a} and \eref{1.34}--\eref{1.34a}, respectively, with the following boundary conditions: 
\begin{equation}\label{1.46}
u^{(1),n}_{0,k}=\kappa_2 u^{(2),n}_{1,k}\quad  {\rm and} \quad  u^{(2),n}_{M+1,k}= \kappa_1 u^{(1),n}_{M,k},
\end{equation}
which are different from \eref{1.32a}--\eref{1.33}.

In order to show the decay of ${\cal L}^n$ and hence the $L^2$-stability, we proceed by combining the proof of Theorem \ref{thm1} and the results in \S \ref{sec4.2}. First, we obtain 
\begin{equation}\label{1.46a}
\begin{aligned}
&  \sum_{j=1}^{M}  \Big( \big(u^{(1),n+1}_{j,k}\big)^2-\big(u^{(1),n}_{j,k}\big)^2  \Big)    e^{-\mu_1 x_j}+ \sum_{j=1}^{M}     \Big( \big(u^{(2),n+1}_{j,k}\big)^2-\big(u^{(2),n}_{j,k}\big)^2  \Big) e^{\mu_2 x_j}   \\
\hspace{-0.3cm}  \le  &\sum_{j=1}^{M} \big(u^{(1),n}_{j,k}\big)^2 \big(D^{(1),n}_{j+1,k}e^{-\mu_1 \dx}-D^{(1),n}_{j,k}\big) e^{-\mu_1 x_j}
 +\sum_{j=1}^{M}  \big(u^{(2),n}_{j,k}\big)^2 \big(D^{(2),n}_{j,k}e^{-\mu_2 \dx}-D^{(2),n}_{j+1,k}\big)e^{\mu_2 x_j} +R,
\end{aligned}
\end{equation}
where 
$$
\begin{aligned}
R&=D^{(1),n}_{1,k} e^{-\mu_1 x_1} \big(u^{(1),n}_{0,k}\big)^2-D^{(1),n}_{M,k}e^{-\mu_1 x_{M+1}}\big(u^{(1),n}_{M,k}\big)^2  
 +D^{(2),n}_{M,k}e^{\mu_2 x_M} \big(u^{(2),n}_{M,k}\big)^2-D^{(2),n}_{1,k}e^{\mu_2x_0}\big(u^{(2),n}_{1,k}\big)^2  \\
 &\overset{\eref{1.46}}{=}\big(u^{(2),n}_{1,k}\big)^2\Big( D^{(1),n}_{1,k}\kappa_2^2 e^{-\mu_1 x_1} -D^{(2),n}_{1,k}e^{\mu_2x_0}   \Big) 
 +\big(u^{(1),n}_{M,k}\big)^2 \Big( D^{(2),n}_{M,k} \kappa_1^2 e^{\mu_2 x_{M}} -D^{(1),n}_{M,k}e^{-\mu_1 x_{M+1}}   \Big) \\
 &\overset{\eref{1.40}}{\le}\big(u^{(2),n}_{1,k}\big)^2\Big( D^{\rm max}_{1}\kappa_2^2 e^{-\mu_1 x_1} -D^{\rm min}_{2}e^{\mu_2x_0}   \Big) 
 +\big(u^{(1),n}_{M,k}\big)^2 \Big( D^{\rm max}_{2}\kappa_1^2 e^{\mu_2 x_{M}}-D^{\rm min}_{1}e^{-\mu_1 x_{M+1}}   \Big) \\
 &=\big(u^{(2),n}_{1,k}\big)^2 D^{\rm min}_2 e^{\mu_2(x_1-\dx)}\Bigg( \bigg( \sqrt{\frac{D^{\rm max}_1}{D^{\rm min}_2}}\kappa_2   \bigg)^2e^{-\mu_1 x_1 -\mu_2x_0} -1 \Bigg) \\
 &+\big(u^{(1),n}_{M,k}\big)^2 D^{\rm min}_{1}e^{-\mu_1 x_{M+1}} \Bigg(  \bigg(\sqrt{\frac{D^{\rm max}_2}{D^{\rm min}_1}} \kappa_1  \bigg)^2 e^{\mu_2 x_{M}+\mu_1 x_{M+1}}-1   \Bigg).
\end{aligned}
$$
For simplicity, we choose $\mu_1 = \mu_2 =\tilde{\mu}$, and hence, provided that 
\begin{equation}\label{1.47}
\kappa_2 < \sqrt{\frac{D^{\rm min}_2}{D^{\rm max}_1}}\quad {\rm and} \quad \kappa_1 < \sqrt{\frac{D^{\rm min}_1}{D^{\rm max}_2}},
\end{equation}
there exists $\tilde{\mu}$,
\begin{equation*}
 0<\tilde{\mu} \le \frac{1}{2+\dx} \Bigg( \ln \Bigg(  \sqrt{\frac{D^{\rm max}_2}{D^{\rm min}_1}}\kappa_1  \Bigg)^{-2}   \Bigg),
\end{equation*}
such that $R \le 0$. The remaining two terms on the RHS of \eref{1.46a} are estimated as in the proof of Theorem \ref{thm1} and in \S \ref{sec4.2} (we omit the details for the sake of brevity), and we finally obtain the results summarized in the next theorem. 
\begin{thm}\label{thm4.3}
Let $m=1$, $p=2$, and the assumptions of  Theorem \ref{thm3} hold. Then, there exist $\mu_1>0$, $\mu_2>0$ such that for the numerical solution $u^{(1),n}_{j,k}$ and $u^{(2),n}_{j,k}$ defined by \eref{1.10} and  \eref{1.34}, respectively,  with the boundary conditions given by \eref{1.46} with $\kappa_1$ and $\kappa_2$ bounded by \eref{1.47},  the Lyapunov function \eref{1.41} satisfies \eref{1.19} for $\nu$ given by \eref{1.40d}. Moreover, $u^{(i),n}_{j,k}$ is exponentially stable in the discrete $L^2$-norm, that is, \eref{1.20} is valid.	
\end{thm}

\begin{rmk}
In this paper, we do not provide a detailed discussion of the well-posedness of random hyperbolic systems. We note that for each fixed
realization of the random variable, the system reduces to a deterministic hyperbolic system of conservation laws with smooth flux functions
and well-posed initial–boundary value problems, for which the classical theory (see, e.g., \cite{DafermosBook}) applies. Since the random
input enters only through the initial and boundary data, the standard well-posedness results can be invoked for each realization, and our
primary focus here is the numerical stabilization analysis based on discrete Lyapunov functions.
\end{rmk}

\section{Numerical Results}\label{sec5}
\setcounter{equation}{0}

In this section, we validate the theoretical estimates in \S \ref{sec3} and \ref{sec4} on the linear advection equation and the linearized shallow-water equations. We also demonstrate the decay rates on an example for the nonlinear shallow-water equations. In all of the examples, we take $K=100$.  In Examples 1--5, we take the CFL number 1, while in Examples 6--8, we take the CFL number 0.5.

\subsection{Boundary Stabilization of the Linear Advection Equation}
We first consider the linear advection equation 
\begin{equation*}
u_t+u_x=0,
\end{equation*}
with $m=p=1$ and $\Lambda_1=1$.

\paragraph{Example 1.} In the first example, we consider the following initial and boundary conditions:
\begin{equation*}
u(x,0,\xi)=-\frac{1}{2}(-\sigma+\xi),   \quad u(0,t,\xi)=0.75 u(1,t,\xi), \quad \xi \in [-\sigma, \sigma],\quad\rho(\xi)=\frac{1}{2\sigma}.  
\end{equation*}

We compute the numerical results until the final time $T=12$ on a sequence of uniform meshes with $\dx=1/100$, 1/200, 1/400, 1/800, and 1/1600, and present the obtained numerical results for $\sigma=\frac{1}{2}$, 1, and 2 in Table \ref{table1}.  One can clearly see that, as expected, the approximate values  $\tilde{\nu}:=-\frac{1}{T}\ln \frac{{\cal L}(T)}{{\cal L}^0}$ converge toward the theoretical values $\nu$ defined by Theorem \ref{thm4.1} as $\dx \to 0$.

\begin{table}[ht!]
\centering
\resizebox{\linewidth}{!}{%
\begin{tabular}{|c|c c c|c c c|c c c|}
\hline
\multirow{2}{*}{$\dx$} & \multicolumn{3}{c|}{$\sigma=\frac{1}{2}$} & \multicolumn{3}{c|}{$\sigma=1$} & \multicolumn{3}{c|}{$\sigma=2$} \\
\cline{2-10}
& $E$ & $\tilde{\nu}$ & $\nu$ & $E$ & $\tilde{\nu}$ & $\nu$ & $E$ & $\tilde{\nu}$ & $\nu$ \\
\hline
$1/100$  & 3.61e-4 & 5.691e-1 & 5.721e-1 & 1.45e-3 & 5.691e-1 & 5.721e-1 & 5.78e-3 & 5.691e-1 & 5.721e-1 \\
$1/200$  & 1.82e-4 & 5.722e-1 & 5.737e-1 & 7.26e-4 & 5.722e-1 & 5.737e-1 & 2.90e-3 & 5.722e-1 & 5.737e-1 \\
$1/400$  & 9.09e-5 & 5.738e-1 & 5.745e-1 & 3.64e-4 & 5.738e-1 & 5.745e-1 & 1.46e-3 & 5.738e-1 & 5.745e-1 \\
$1/800$  & 4.55e-5 & 5.746e-1 & 5.750e-1 & 1.82e-4 & 5.746e-1 & 5.750e-1 & 7.28e-4 & 5.746e-1 & 5.750e-1 \\
$1/1600$ & 2.28e-5 & 5.750e-1 & 5.752e-1 & 9.11e-5 & 5.750e-1 & 5.752e-1 & 3.64e-4 & 5.750e-1 & 5.752e-1 \\
\hline
\end{tabular}
}
\caption{\sf Example 1: The decay of $E:=||e^{-\nu t}   {\cal L}^{0}-{\cal L}||_\infty$ together with the corresponding values of $\tilde{\nu}$ and $\nu$.}
\label{table1}
\end{table}

\paragraph{Example 2.} In the second example, we use the same set-up as in Example 1, but with the following discontinuous initial data:
\begin{equation*}
u(x,0,\xi)=\begin{cases}
	        -\frac{1}{2}, & \text{$x<\frac{1}{4}$},\\[1.ex]
	        -\frac{1}{2}(-\sigma+\xi), & \text{otherwise}.
\end{cases}
\end{equation*}

We also compute the numerical results until the final time $T=12$ on the same sequence of uniform meshes with $\dx=1/100$, 1/200, 1/400, 1/800, and 1/1600.  The obtained numerical results for $\sigma=\frac{1}{2}$, 1, and 2 are presented in Table \ref{table2}. As in Example 1, one can clearly see that the approximate values  $\tilde{\nu}$ converge toward the theoretical ones as $\dx \to 0$. 

\begin{table}[ht!]
\centering
\resizebox{\linewidth}{!}{%
\begin{tabular}{|c|c c c|c c c|c c c|}
\hline
\multirow{2}{*}{$\dx$} & \multicolumn{3}{c|}{$\sigma=\frac{1}{2}$} & \multicolumn{3}{c|}{$\sigma=1$} & \multicolumn{3}{c|}{$\sigma=2$} \\
\cline{2-10}
& $E$ & $\tilde{\nu}$ & $\nu$ & $E$ & $\tilde{\nu}$ & $\nu$ & $E$ & $\tilde{\nu}$ & $\nu$ \\
\hline
$1/100$  & 2.37e-3 & 5.691e-1 & 5.721e-1 & 2.02e-3 & 5.691e-1 & 5.721e-1 & 6.62e-4 & 5.691e-1 & 5.721e-1 \\
$1/200$  & 1.19e-3 & 5.722e-1 & 5.737e-1 & 1.01e-3 & 5.722e-1 & 5.737e-1 & 3.34e-4 & 5.722e-1 & 5.737e-1 \\
$1/400$  & 5.95e-4 & 5.738e-1 & 5.745e-1 & 5.07e-4 & 5.738e-1 & 5.745e-1 & 1.68e-4 & 5.738e-1 & 5.745e-1 \\
$1/800$  & 2.98e-4 & 5.746e-1 & 5.750e-1 & 2.54e-4 & 5.746e-1 & 5.750e-1 & 8.40e-5 & 5.746e-1 & 5.750e-1 \\
$1/1600$ & 1.49e-4 & 5.750e-1 & 5.752e-1 & 1.27e-4 & 5.750e-1 & 5.752e-1 & 4.22e-5 & 5.750e-1 & 5.752e-1 \\
\hline
\end{tabular}
}
\caption{\sf Example 2: The decay of $E=||e^{-\nu t}   {\cal L}^{0}-{\cal L}||_\infty$  together with the corresponding values of $\tilde{\nu}$ and $\nu$.}
\label{table2}
\end{table}

\subsection{Boundary Stabilization of the Linearized Saint-Venant System}\label{sec5.2}
In this section, we consider the boundary damping for the linearized Saint-Venant system of shallow-water equations, which is inspired by many recent results on the continuous formulation; see, e.g., \cite{Coron1999,bastin2007lyapunov,coron2007strict,deHalleux2003}. The primary goal is to control the water depth and velocity in an open canal or a network of canals. External factors such as weather variations or uncontrolled inflows can cause fluctuations in water depth, which must be mitigated to maintain a predetermined target depth $\xbar h$ and velocity $\xbar v$. As in \cite{coron2007strict,deHalleux2003}, we neglect the source terms. The target water state is a constant with $\xbar h(x)=4$ and $\xbar v(x)=\frac{5}{2}$. To analyze small deviations from this steady state, the Saint-Venant system of shallow-water equations
\begin{equation*}
 \begin{pmatrix}
	 h \\
	 q
 \end{pmatrix} _t+ \begin{pmatrix}
	q  \\
	\dfrac{q^2}{h}+\dfrac{g}{2}h^2
 \end{pmatrix}_x=\begin{pmatrix}
	0 \\
	0 
 \end{pmatrix},
\end{equation*}
where $q=hv$, is linearized by introducing perturbations $(\delta h, \delta v)$ with $h=\xbar h+\delta h$ and $v = \xbar v +\delta v$. The resulting diagonalized equations for these perturbations are given by:
\begin{equation}\label{4.5}
 \begin{pmatrix}
	u^{(1)} \\
	u^{(2)}
 \end{pmatrix} _t+ \begin{pmatrix}
	\Lambda_1 & 0  \\
	0 & \Lambda_2
 \end{pmatrix}\begin{pmatrix}
	u^{(1)}\\
	u^{(2)} 
 \end{pmatrix} _x=\begin{pmatrix}
	0 \\
	0 
 \end{pmatrix},
\end{equation}
where  
\begin{equation}\label{4.5a}
\Lambda_1=\xbar v + \sqrt{\mathstrut{g\,\xbar h}}, \quad \Lambda_2=\xbar v - \sqrt{\mathstrut{g\,\xbar h}}, \quad   u^{(1)}=\delta v+\sqrt{\frac{g}{\xbar h}}\delta h, \quad  {\rm and} \quad  u^{(2)}=\delta v-\sqrt{\frac{g}{\xbar h}}\delta h.
\end{equation}

\paragraph{Example 3.} In this example, we consider the diagonalized equations  \eref{4.5}--\eref{4.5a} with the following initial conditions:
\begin{equation*}
	\delta h(x,0,\xi)=\frac{1}{2}\sin(\pi x)(-\sigma+\xi), \quad \delta v(x,0,\xi)=\frac{20}{8+\sin(\pi x)}-\frac{5}{2}, \quad \xi \in [-\sigma, \sigma],\quad \rho(\xi)=\frac{1}{2\sigma},
\end{equation*}
subject to the boundary conditions $u^{(1)}(x,0,\xi) = 0.8\, u^{(1)}(x,1,\xi)$ and $u^{(2)}(x,1,\xi)  =0.8\, u^{(2)}(x,0,\xi).$

We compute the numerical results until the final time $T=6$ on  a sequence of uniform meshes with $\dx=1/100$, 1/200, 1/400, 1/800, and 1/1600, and present the obtained numerical results for $\sigma=\frac{1}{2}$, 1, and 2 in Table \ref{table3}.  The numerical experiments validate the theoretical results, confirming that the applied feedback boundary control effectively stabilizes the system.

\begin{table}[ht!]
\centering
%\resizebox{\linewidth}{!}{%
\begin{tabular}{|c|ccc|ccc|ccc|}
\hline
\multirow{2}{1em}{$\dx$}&\multicolumn{3}{c|}{$\sigma=\frac{1}{2}$}&\multicolumn{3}{c|}{$\sigma=1$}&\multicolumn{3}{c|}{$\sigma=2$}\\
\cline{2-10}&$E$&$\tilde{\nu}$&$\nu$&$E$&$\tilde{\nu}$&$\nu$&$E$&$\tilde{\nu}$&$\nu$\\
\hline
$1/100$  &1.23e-2&1.806&1.699& 2.69e-2 &1.807&1.699&8.53e-2 &1.808&1.699\\
$1/200$  &1.20e-2&1.810&1.703& 2.63e-2 &1.811&1.703&8.33e-2 &1.812&1.703\\
$1/400$  &1.18e-2&1.807&1.705& 2.59e-2 &1.808&1.705&8.21e-2 &1.809&1.705\\
$1/800$  &1.17e-2&1.801&1.706& 2.57e-2 &1.801&1.706&8.14e-2 &1.802&1.706\\
$1/1600$ &1.16e-2&1.796&1.706& 2.55e-2 &1.796&1.706&8.10e-2 &1.797&1.706\\
\hline
\end{tabular}
%}
\caption{\sf Example 3: The decay of $E=||e^{-\nu t}   {\cal L}^{0}-{\cal L}||_\infty$  together with the corresponding values of $\tilde{\nu}$ and $\nu$.\label{table3}}
\end{table}

\paragraph{Example 4.} In this example, we use the same set-up as in Example 3, but with a different, random velocity perturbation:
\begin{equation*}
\delta v(x,0,\xi)=\frac{10}{4+\delta h}-\frac{5}{2}.
\end{equation*}

We compute the numerical results until the final time $T=6$ on the same sequence of uniform meshes with $\dx=1/100$, 1/200, 1/400, 1/800, and 1/1600, and present the obtained numerical results for $\sigma=\frac{1}{2}$, 1, and 2 in Table \ref{table4}.  The results support the theoretical predictions, demonstrating that the applied boundary control successfully stabilizes the system.

\begin{table}[ht!]
\centering
%\resizebox{\linewidth}{!}{%
\begin{tabular}{|c|ccc|ccc|ccc|}
\hline
\multirow{2}{1em}{$\dx$}&\multicolumn{3}{c|}{$\sigma=\frac{1}{2}$}&\multicolumn{3}{c|}{$\sigma=1$}&\multicolumn{3}{c|}{$\sigma=2$}\\
\cline{2-10}&$E$&$\tilde{\nu}$&$\nu$&$E$&$\tilde{\nu}$&$\nu$&$E$&$\tilde{\nu}$&$\nu$\\
\hline
$1/100$  &5.36e-3&1.746&1.699& 2.15e-2 &1.746&1.699&8.67e-2 &1.746&1.699\\
$1/200$  &5.15e-3&1.750&1.703& 2.06e-2 &1.750&1.703&8.32e-2 &1.750&1.703\\
$1/400$  &5.02e-3&1.747&1.705& 2.01e-2 &1.747&1.705&8.11e-2 &1.746&1.705\\
$1/800$  &4.95e-3&1.740&1.706& 1.98e-2 &1.740&1.706&7.99e-2 &1.740&1.706\\
$1/1600$ &4.91e-3&1.734&1.706& 1.97e-2 &1.734&1.706&7.93e-2 &1.734&1.706\\
\hline
\end{tabular}
%}
\caption{\sf Example 4: The decay of $E=||e^{-\nu t}   {\cal L}^{0}-{\cal L}||_\infty$ together with the corresponding values of $\tilde{\nu}$ and $\nu$.\label{table4}}
\end{table}

\paragraph{Example 5.} In this example, we use the same set-up as in Example 4, but with different boundary conditions:
\begin{equation*}
u^{(1)}(x,0,\xi) = 0.6 u^{(2)}(x,0,\xi), \quad  u^{(2)}(x,1,\xi) = 0.6u^{(1)}(x,1,\xi).
\end{equation*}

We compute the numerical results until the final time $T=6$ on a sequence of uniform meshes with $\dx=1/100$, 1/200, 1/400, 1/800, and 1/1600 and present the obtained numerical results for $\sigma=\frac{1}{2}$, 1, and 2 in Table \ref{table2.4}.  As in the previous examples, the applied boundary control stabilized the system. We also stress that this example validates the results stated in Theorem \ref{thm4.3}. 

\begin{table}[ht!]
\centering
%\resizebox{\linewidth}{!}{%
\begin{tabular}{|c|ccc|ccc|ccc|}
\hline
\multirow{2}{1em}{$\dx$}&\multicolumn{3}{c|}{$\sigma=\frac{1}{2}$}&\multicolumn{3}{c|}{$\sigma=1$}&\multicolumn{3}{c|}{$\sigma=2$}\\
\cline{2-10}&$E$&$\tilde{\nu}$&$\nu$&$E$&$\tilde{\nu}$&$\nu$&$E$&$\tilde{\nu}$&$\nu$\\
\hline
$1/100$  &9.38e-3&5.536&1.751& 3.76e-2 &5.536&1.751&1.52e-1 &5.535&1.751\\
$1/200$  &9.24e-3&5.511&1.764& 3.70e-2 &5.510&1.764&1.50e-1 &5.510&1.764\\
$1/400$  &9.17e-3&5.485&1.770& 3.67e-2 &5.485&1.770&1.49e-1 &5.484&1.770\\
$1/800$  &9.13e-3&5.469&1.773& 3.66e-2 &5.469&1.773&1.48e-1 &5.469&1.773\\
$1/1600$ &9.11e-3&5.460&1.775& 3.65e-2 &5.460&1.775&1.48e-1 &5.460&1.775\\
\hline
\end{tabular}
%}
\caption{\sf Example 5: The decay of $E=||e^{-\nu t}   {\cal L}^{0}-{\cal L}||_\infty$ together with the corresponding values of $\tilde{\nu}$ and $\nu$.\label{table2.4}}
\end{table}

\subsection{Boundary Stabilization of the Nonlinear Equations}
In this section, we consider the same setting as in \S \ref{sec5.2}, but with different equations for the perturbations:
\begin{equation}\label{5.5aa}
 \begin{pmatrix}
	u^{(1)} \\
	u^{(2)}
 \end{pmatrix} _t+ \begin{pmatrix}
	\xbar v +\delta v + \sqrt{g \big(\xbar h+\delta h\big)} & 0  \\
	0 & \xbar v+\delta v - \sqrt{g \big(\xbar h+\delta h\big)}
 \end{pmatrix}\begin{pmatrix}
	u^{(1)}\\
	u^{(2)}
 \end{pmatrix} _x=\begin{pmatrix}
	0 \\
	0 
 \end{pmatrix},
\end{equation}
where $u^{(1)}$ and $u^{(2)}$ are defined by \eref{4.5a}. Notice that the system \eref{5.5aa} is nonlinear. 

\paragraph{Example 6.} In this example, we take the same setting as in Example 4 and compute the numerical results until the final time
$T=6$ on a sequence of uniform meshes with $\dx=1/100$, $1/200$, $1/400$, $1/800$, and $1/1600$. The obtained results, presented in Table
\ref{table6.1}, demonstrate that the proposed boundary control is capable of stabilizing the nonlinear system.  We stress,
however, that even though the numerical results demonstrate stability, a rigorous proof establishing the convergence properties of the CU
scheme remains an open question. Note that the exponential stability of a general second-order slope limiter scheme for scalar conservation
laws subject to a dissipative boundary condition has been recently studied in \cite{Dus2022}.
\begin{table}[ht!]
\centering
%\resizebox{\linewidth}{!}{%
\begin{tabular}{|c|cc|cc|cc|ccc|}
\hline
\multirow{2}{1em}{$\dx$}&\multicolumn{2}{c|}{$\sigma=\frac{1}{2}$}&\multicolumn{2}{c|}{$\sigma=1$}&\multicolumn{2}{c|}{$\sigma=2$}\\
\cline{2-7}&$E$&$\tilde{\nu}$&$E$&$\tilde{\nu}$&$E$&$\tilde{\nu}$\\
\hline
$1/100$  &6.92e-3&3.910 &  3.91e-2 &3.907&1.09e-1 &3.892\\
$1/200$  &6.71e-3&3.931 &  2.68e-2 &3.928&1.06e-1 &3.912\\
$1/400$  &6.59e-3&3.940 &  2.63e-2 &3.935&1.04e-1 &3.919\\
$1/800$  &6.52e-3&3.935 &  2.60e-2 &3.932&1.03e-1 &3.919\\
$1/1600$ &6.49e-3&3.930 &  2.59e-2 &3.928&1.02e-1 &3.915\\
\hline
\end{tabular}
%}
\caption{\sf Example 6: The decay of $E=||e^{-\nu t}{\cal L}^0-{\cal L}||_\infty$ together with the corresponding values of $\tilde{\nu}$.
\label{table6.1}}
\end{table}

\subsection{Second-Order Discretization of the Linearized Saint-Venant System }
In this section, we apply the second-order semi-discrete central-upwind (CU) scheme from \cite{Kurganov07,KNP} to the linearized Saint-Venant system \eref{4.5}--\eref{4.5a} to demonstrate that faster convergence can be achieved when higher-order schemes are used. A brief overview of the 1-D CU scheme is provided in Appendix \ref{appa}, while details on the computation of boundary conditions are discussed in Appendix \ref{appb}.

\paragraph{Example 7.} In this example, we take the same setting as in Example 4 and compute the numerical results using the CU scheme until the final time $T=6$ on uniform meshes with $\dx=1/100$, 1/200, 1/400, 1/800, and 1/1600. We present the obtained numerical results for $\sigma=\frac{1}{2}$, 1, and 2 in Table \ref{table7}. Compared with the corresponding first-order results (see Table \ref{table4}), the decay of $E$ is now much faster and the values of $\tilde{\nu}$ are closer to the theoretical estimates $\nu$. However, a rigorous proof of discrete convergence for the CU scheme or any other higher-order schemes remains an open problem.  
\begin{table}[ht!]
\centering
%\resizebox{\linewidth}{!}{%
\begin{tabular}{|c|ccc|ccc|ccc|}
\hline
\multirow{2}{1em}{$\dx$}&\multicolumn{3}{c|}{$\sigma=\frac{1}{2}$}&\multicolumn{3}{c|}{$\sigma=1$}&\multicolumn{3}{c|}{$\sigma=2$}\\
\cline{2-10}&$E$&$\tilde{\nu}$&$\nu$&$E$&$\tilde{\nu}$&$\nu$&$E$&$\tilde{\nu}$&$\nu$\\
\hline
$1/100$  &1.08e-4&1.702&1.699& 4.33e-4 &1.702&1.699&1.75e-3 &1.704&1.699\\
$1/200$  &1.07e-4&1.704&1.703& 4.29e-4 &1.704&1.703&1.73e-3 &1.705&1.703\\
$1/400$  &1.06e-4&1.712&1.705& 4.26e-4 &1.712&1.705&1.72e-3 &1.712&1.705\\
$1/800$  &1.06e-4&1.718&1.706& 4.25e-4 &1.718&1.706&1.72e-3 &1.718&1.706\\
$1/1600$ &1.06e-4&1.718&1.706& 4.25e-4 &1.718&1.706&1.72e-3 &1.718&1.706\\
\hline
\end{tabular}
%}
\caption{\sf Example 7: The decay of $E=||e^{-\nu t}   {\cal L}^{0}-{\cal L}||_\infty$ together with the corresponding values of $\tilde{\nu}$ and $\nu$.\label{table7}}
\end{table}

\subsection{Boundary Stabilization of the Linearized Saint-Venant System with Source Terms}
In this section, we consider a numerical example of the boundary damping for the linearized Saint-Venant system with source terms. We take
the same setting as in \S\ref{sec5.2}, but with different equations, which read as
\begin{equation*}
\begin{pmatrix}u^{(1)}\\u^{(2)}\end{pmatrix}_t+\begin{pmatrix}\Lambda_1&0\\0&\Lambda_2\end{pmatrix}
\begin{pmatrix}u^{(1)}\\u^{(2)}\end{pmatrix}_x=-\begin{pmatrix}0.1&0\\0&0.1\end{pmatrix}\begin{pmatrix}u^{(1)}\\u^{(2)}\end{pmatrix},
\end{equation*}
where $u^{(1)}$ and $u^{(2)}$ are defined by \eref{4.5a}.

\paragraph{Example 8.} In this example, we use the same set-up as in Example 4, compute the numerical results until the final time $T=6$ on
the same sequence of uniform meshes with $\dx=1/100$, $1/200$, $1/400$, $1/800$, and $1/1600$, and present the obtained numerical results
for $\sigma=\hf$, $1$, and $2$ in Table \ref{table8}. The results support the theoretical predictions provided in Appendix \ref{appc},
demonstrating that the applied boundary control successfully stabilizes the system.
\begin{table}[ht!]
\centering
\begin{tabular}{|c|ccc|ccc|ccc|}
\hline
\multirow{2}{1em}{$\dx$}&\multicolumn{3}{c|}{$\sigma=\frac{1}{2}$}&\multicolumn{3}{c|}{$\sigma=1$}&\multicolumn{3}{c|}{$\sigma=2$}\\
\cline{2-10}&$E$&$\tilde{\nu}$&$\nu$&$E$&$\tilde{\nu}$&$\nu$&$E$&$\tilde{\nu}$&$\nu$\\
\hline
$1/100$ &1.09e-4&1.946&1.699&4.36e-4&1.946&1.699&1.76e-3&1.946&1.699\\
$1/200$ &1.07e-4&1.952&1.703&4.31e-4&1.952&1.703&1.74e-3&1.952&1.703\\
$1/400$ &1.07e-4&1.951&1.705&4.28e-4&1.951&1.705&1.73e-3&1.950&1.705\\
$1/800$ &1.07e-4&1.944&1.706&4.27e-4&1.944&1.706&1.72e-3&1.944&1.706\\
$1/1600$&1.06e-4&1.937&1.706&4.26e-4&1.937&1.706&1.72e-3&1.937&1.706\\
\hline
\end{tabular}
\caption{\sf Example 8: The decay of $E=||e^{-\nu t}{\cal L}^0-{\cal L}||_\infty$ together with the corresponding values of $\tilde\nu$ and
$\nu$.\label{table8}}
\end{table}

\section{Conclusion}\label{sec6}
In this paper, we have extended the Lyapunov-based stabilization framework to random systems of hyperbolic conservation laws,
where uncertainties arise from boundary controls and initial data. By integrating a stochastic discrete Lyapunov function into the numerical
analysis, we rigorously established exponential stability for the finite-volume discretization of these systems. Theoretical decay rates
were derived, demonstrating their dependence on boundary control parameters, grid resolution, and the statistical properties of random
inputs. We experimentally checked the stability of one nonlinear system, and the obtained results demonstrate convergence as well even
though the rigorous proof is currently out of reach. We have applied the second-order central-upwind (CU) scheme to the linearized
Saint-Venant system of shallow-water equations and compared its performance with the performance of the first-order upwind scheme. The
conducted numerical experiments confirm that the CU scheme achieves much smaller values of $||e^{-\nu t}{\cal L}^0-{\cal L}||_\infty$ (which
means that the obtained stability estimates are sharper) while maintaining stability. We have also experimentally checked the
stability in one numerical example of the random hyperbolic system of balance laws, which demonstrates that the applied boundary control
successfully stabilizes the system. Future research directions include extending this framework to more complex nonlinear (nonconservative)
systems and high-order numerical schemes.

\begin{DA}
\paragraph{Funding.} The work of S. Chu and  M. Herty was funded by the Deutsche Forschungsgemeinschaft (DFG, German Research Foundation) - SPP 2410 Hyperbolic Balance Laws in Fluid Mechanics: Complexity, Scales, Randomness (CoScaRa) within the Project(s) HE5386/26-1 (Numerische Verfahren für gekoppelte Mehrskalenprobleme,525842915) and (Zufällige kompressible Euler Gleichungen: Numerik und ihre Analysis, 525853336) HE5386/27-1, and  the Deutsche Forschungsgemeinschaft (DFG, German Research Foundation) - SPP 2183: Eigenschaftsgeregelte Umformprozesse with the Project(s) HE5386/19-2,19-3 Entwicklung eines flexiblen isothermen Reckschmiedeprozesses für die eigenschaftsgeregelte Herstellung von Turbinenschaufeln aus Hochtemperaturwerkstoffen (424334423). The work of A. Kurganov was supported in part by NSFC grants 12171226 and W2431004.

\paragraph{Conflicts of interest.} On behalf of all authors, the corresponding author states that there is no conflict of interest.

\paragraph{Data and software availability.} The data that support the findings of this study and FORTRAN codes developed by the authors and
used to obtain all of the presented numerical results are available from the corresponding author upon reasonable request.
\end{DA}

\appendix
\setcounter{equation}{0}
\section{Central-Upwind (CU) Scheme}\label{appa}

In this appendix, we briefly describe the semi-discrete CU scheme, which was introduced in \cite{Kurganov00,KNP} (see also \cite{CCHKL_22,CKX_24,KX_22} for recent low-dissipation modifications of the CU schemes)  as a ``black-box'' solver for general hyperbolic systems of conservation laws.

As CU schemes are finite-volume methods, the computed quantities are the cell averages, $\xbar \bmu_{\jph,k}(t)\approx\frac{1}{\dx}\int_{x_j}^{x_{j+1}} \bmu(x,t) \mathrm{d}x$, $j=0,\ldots,M-1$. The semi-discrete CU scheme from \cite{KNP} for the system \eref{1.1} reads as
\begin{equation}\label{a1}
\frac{\mathrm{d} \overline{\bmu}_{\jph,k}(t)}{\mathrm{d}t}=-\frac{{{\bmF}}_{j+1,k}(t)-{{\bmF}}_{j,k}(t)}{\Delta x},
\end{equation}
where the CU numerical fluxes are given by 
\begin{equation}\label{a2}
{\bmF}_{j,k}(t)=\frac{a^{+}_{j,k}(t) \mf\big(\bmu^{-}_{j,k}(t)\big)
-a^{-}_{j,k}(t) \mf\big(\bmu^{+}_{j,k}(t)\big)}{a^{+}_{j,k}(t)-a^{-}_{j,k}(t)}
+\frac{a^{+}_{j,k}(t) a^{-}_{j,k}(t)}{a^{+}_{j,k}(t)
-a^{-}_{j,k}(t)}\big(\bmu^{+}_{j,k}(t)-\bmu^{-}_{j,k}(t)\big).
\end{equation}
Here, $\bmu^{\pm}_{j,k}(t)$ are the right/left-sided values of a global (in space) piecewise linear interpolant
\begin{equation}
\widetilde{\bmu}(x;t)=  \widetilde \bmu_{\jph,k}(x;t):=\,\xbar \bmu_{\jph,k}(t)+\big(\bmu_x(t)\big)_{\jph,k}(x-x_\jph),\quad x\in (x_j,x_{j+1}),
\label{2.3}
\end{equation}
at the cell interface $x=x_j$, namely,
\begin{equation*}
\bmu^{-}_{j,k}(t)=\,\xbar \bmu_{\jmh,k}(t)+\frac{\dx}{2}\big(\bmu_x(t)\big)_{\jmh,k},\quad \bmu^{+}_{j,k}(t)=\,\xbar \bmu_{\jph,k}(t)-\frac{\dx}{2}\big(\bmu_x(t)\big)_{\jph,k}.
\end{equation*}
In order to ensure a non-oscillatory nature of this reconstruction, one needs to compute the slopes $\big(\bmu_x(t)\big)_{\jph,k}$ in \eref{2.3} using a nonlinear limiter. In the numerical experiment reported in \S\ref{sec4}, we have used a minmod limiter
\cite{lie03,Sweby84}:
\begin{equation}\label{a4}
\big(\bmu_x(t)\big)_{\jph,k}={\rm minmod}\bigg(\frac{\,\xbar \bmu_{\jph,k}(t)-\,\xbar \bmu_{\jmh,k}(t)}{\dx},\,\,\frac{\,\xbar \bmu_{j+\frac{3}{2},k}(t)-\,\xbar \bmu_{\jph,k}(t)}{\dx}\bigg), \quad j=1,\ldots,M-2. 
\end{equation}
Here, the minmod function is defined as
\begin{equation*}
{\rm minmod}(z_1,z_2,\ldots):=\begin{cases}
\min_j\{z_i\}&\mbox{if}~z_i>0\quad\forall\,i,\\
\max_j\{z_i\}&\mbox{if}~z_i<0\quad\forall\,i,\\
0            &\mbox{otherwise},
\end{cases}
\end{equation*}
and applied in a component-wise manner.

We note that in order to apply the minmod limiter \eref{a4} at the boundary cells, one needs to use the ghost values $\xbar \bmu_{-\frac{1}{2},k}(t)$ and  $\xbar \bmu_{M+\frac{1}{2},k}(t)$, which may be unavailable for certain boundary conditions including those considered in this paper. Therefore, at the boundary cells, one can use a different minmod limiting:
\begin{equation}\label{a4a}
\begin{aligned}
&\big(\bmu_x(t)\big)_{\frac{1}{2},k}={\rm minmod}\Bigg(\frac{\xbar \bmu_{\frac{1}{2},k}(t)-\bmu^{-}_{0,k}(t)}{\dx/2},\,\,\frac{\xbar \bmu_{\frac{3}{2},k}(t)-\xbar \bmu_{\frac{1}{2},k}(t)}{\dx}  \Bigg),\\
&\big(\bmu_x(t)\big)_{M-\frac{1}{2},k}={\rm minmod}\Bigg( \frac{\xbar \bmu_{M-\frac{1}{2},k}(t)-\xbar \bmu_{M-\frac{3}{2},k}(t)}{\dx},\,\, \frac{\bmu^{+}_{M,k}(t)-\xbar \bmu_{M-\frac{1}{2},k}(t)}{\dx/2}  \Bigg),
\end{aligned}
\end{equation}
where $\bmu^{-}_{0,k}(t)$ and $\bmu^{+}_{M,k}(t)$ are the one-sided boundary point values, which can be computed using the prescribed boundary conditions.

The one-sided local speeds of propagation $a^\pm_\jph(t)$ are estimated using the largest and the smallest eigenvalues of the Jacobian $F(u): \Lambda_p <  \ldots< \Lambda_{1}$. This can be done, for example, by taking
\begin{equation*}
a^{+}_{j,k}(t)=\max\big\{\Lambda_1\big(\bmu^{+}_{j,k}(t)\big),\Lambda_1\big(\bmu^{-}_{j,k}(t)\big),0\big\},  \quad 
a^{-}_{j,k}(t)=\min\big\{\Lambda_p\big(\bmu^{+}_{j,k}(t)\big),\Lambda_p\big(\bmu^{-}_{j,k}(t)\big),0\big\}.
\end{equation*}

Finally, the ODE system \eref{a1} has to be numerically integrated using a stable and sufficiently accurate ODE solver. In the numerical experiment reported in \S \ref{sec4}, we have used the three-stage third-order strong stability preserving Runge-Kutta (SSP-RK3) method (see, e.g., \cite{Gottlieb11,Gottlieb12}) and use the CFL number 0.45. 

\section{Boundary Conditions for the CU Scheme}\label{appb}
In this appendix, we consider the linearized  Saint-Venant system of shallow-water equations \eref{4.5}--\eref{4.5a} to show how to compute the boundary values $\bmu^\pm_{M,k}$ and $\bmu^\pm_{0,k}$ needed to evaluate the CU numerical fluxes \eref{a2} at the boundaries. 

We assume that the solution (its cell averages and all of the cell interface values, including the boundary ones) is available at all times until a certain time level $t$. We also assume that at the time level $t+\tau$ (for the SSP-RK3 ODE solver $\tau$ may be either $\frac{\dt}{2}$ or $\dt$) the cell averages $\xbar \bmu_{\jph,k}(t+\tau)$ have been already computed and our goal now is to obtain $u^{(1),\,\pm}_{M,k}(t+\tau)$ and  $u^{(2),\,\pm}_{0,k}(t+\tau)$ since $u^{(1),\,\pm}_{0,k}(t+\tau)$ and $u^{(2),\,\pm}_{M,k}(t+\tau)$ can then be directly computed from the prescribed boundary conditions. Here, we show how to compute $u^{(1),\,\pm}_{M,k}(t+\tau)$ since $u^{(2),\,\pm}_{0,k}(t+\tau)$ can be obtained in a similar manner. 

We begin with $u^{(1),\,+}_{M,k}(t+\tau)$, which can be calculated using the method of characteristics, namely, 
\begin{equation}\label{b1}
u^{(1),\,+}_{M,k}(t+\tau) =\widetilde{u}^{\,(1)}(x_M-\tau \Lambda_1,t).
\end{equation}
Due to the finite speed of propagation and the CFL restriction on $\tau$, the RHS of \eref{b1} can be easily computed, which results in 
\begin{equation*}
u^{(1),\,+}_{M,k}(t+\tau)= \xbar u^{\,(1)}_{M-\frac{1}{2},k}(t)+\big(u^{\,(1)}_x(t)\big)_{M-\frac{1}{2},k}\bigg(\frac{\dx}{2}- \tau\Lambda_1\bigg).
\end{equation*}

To obtain the left-sided value $u^{(1),\,-}_{M,k}(t+\tau)$, we reconstruct the linear piece in the cell $(x_{M-1},x_M)$:
\begin{equation}\label{b1a}
\widetilde{u}^{\,(1),\,-}_{M-\frac{1}{2},k}(t+\tau) = \xbar u^{\,(1)}_{M-\frac{1}{2},k}(t+\tau)+  (u^{\,(1)}_x(t+\tau))_{M-\frac{1}{2},k}(x-x_{M-\frac{1}{2}}), 
\end{equation}
where the slope is computed as in \eref{a4a}, namely, by
{\small
\begin{equation}\label{b2}
 \big(u^{\,(1)}_x(t+\tau)\big)_{M-\frac{1}{2},k}={\rm minmod}\Bigg( \frac{\xbar u^{\,(1)}_{M-\frac{1}{2},k}(t+\tau)-\xbar u^{\,(1)}_{M-\frac{3}{2},k}(t+\tau)}{\dx},\,\, \frac{u^{(1),+}_{M,k}(t+\tau)-\xbar u^{\,(1)}_{M-\frac{1}{2},k}(t+\tau) }{\dx/2}  \Bigg).
\end{equation}
}
Finally, substituting $x=x_M$ into \eref{b1a}--\eref{b2} results in
\begin{equation*}
u^{(1),\,-}_{M,k}(t+\tau) = \xbar u^{\,(1)}_{M-\frac{1}{2},k}(t+\tau)+ \frac{\dx}{2} \big(u^{\,(1)}_x(t+\tau)\big)_{M-\frac{1}{2},k}.
\end{equation*}

\section{Extension of Theorem \ref{thm3} to Random Hyperbolic Systems of Balance Laws}\label{appc}
In this appendix, we study the $L^2$-stabilization of 1-D random hyperbolic systems of balance laws. We consider a simple case, which is
given by  
\begin{equation}
\bmu_t+F(\bmu)\bmu_x=-B\bmu,
\label{c1}
\end{equation}
where the matrix $B$ is diagonally positive definite, that is,
\begin{equation}
B={\rm diag}(b^{(1)},\dots,b^{(m)}),\qquad b^{(i)}>0.
\label{c1a}
\end{equation}
Making the same assumptions as in Theorem \ref{thm3}, the corresponding upwind scheme for the system \eref{c1}--\eref{c1a} reads as
\begin{equation}
\begin{aligned}
u_{j,k}^{(i),n+1}&=u_{j,k}^{(i),n}-\frac{\dt}{\dx}\,\Lambda_i\!\left(\bmu_{j,k}^n\right)\Big(u_{j,k}^{(i),n}-u_{j-1,k}^{(i),n}\Big)-
\dt\,b^{(i)}u_{j,k}^{(i),n},\\
u_{0,k}^{(i),n+1}&=\kappa_i\,u_{M,k}^{(i),n+1},\qquad i=1,\dots,m,
\end{aligned}
\label{c2}
\end{equation}
and  
\begin{equation}
\begin{aligned}
u_{j,k}^{(i),n+1}&=u_{j,k}^{(i),n}-\frac{\dt}{\dx}\,\Lambda_i\!\left(\bmu_{j+1,k}^n\right)\Big(u_{j+1,k}^{(i),n}-u_{j,k}^{(i),n}\Big)-
\dt\,b^{(i)}u_{j,k}^{(i),n},\\
u_{M+1,k}^{(i),n+1}&=\kappa_i\,u_{1,k}^{(i),n+1},\qquad i=m+1,\dots,p.
\end{aligned}  
\label{c3}
\end{equation}
The scheme \eref{c2}--\eref{c3} can be equivalently rewritten as  
\begin{equation}
\begin{aligned}
u^{(i),n+1}_{j,k}&=u^{(i),n}_{j,k}\Big(1-D^{(i),n}_{j,k}-\dt\,b^{(i)}\Big)+u^{(i),n}_{j-1,k}D^{(i),n}_{j,k},\\ 
u^{(i),n+1}_{0,k}&=\kappa_i\,u^{(i),n+1}_{M,k},\qquad i=1,\dots,m,
\end{aligned}
\label{c4}
\end{equation}  
and
\begin{equation}
\begin{aligned}
u^{(i),n+1}_{j,k}&=u^{(i),n}_{j,k}\Big(1-D^{(i),n}_{j+1,k}-\dt\,b^{(i)}\Big)+u^{(i),n}_{j+1,k}\,D^{(i),n}_{j+1,k},\\
u^{(i),n+1}_{M+1,k}&=\kappa_i\,u^{(i),n+1}_{1,k},\qquad i=m+1,\dots,p.
\end{aligned}
\label{c5}
\end{equation}  
Applying Jensen’s inequality to \eref{c4}--\eref{c5} leads to
\begin{equation*}
\begin{aligned}
\big(u^{(i),n+1}_{j,k}\big)^2&\le(1-\dt\,b^{(i)})\Big[\big(1-D^{(i),n}_{j,k}-\dt\,b^{(i)}\big)\big(u^{(i),n}_{j,k}\big)^2+
D^{(i),n}_{j,k}\big(u^{(i),n}_{j-1,k}\big)^2\Big]\\
&\le\big(1-D^{(i),n}_{j,k}\big)\big(u^{(i),n}_{j,k}\big)^2+D^{(i),n}_{j,k}\big(u^{(i),n}_{j-1,k}\big)^2,\qquad i=1,\dots,m,
\end{aligned}
\end{equation*}
and  
\begin{equation*}
\begin{aligned}
\big(u^{(i),n+1}_{j,k}\big)^2&\le(1-\dt\,b^{(i)})\Big[\big(1-D^{(i),n}_{j+1,k}-\dt\,b^{(i)}\big)\big(u^{(i),n}_{j,k}\big)^2+
D^{(i),n}_{j+1,k}\big(u^{(i),n}_{j+1,k}\big)^2 \Big]\\
&\le\big(1-D^{(i),n}_{j+1,k}\big)\big(u^{(i),n}_{j,k}\big)^2+D^{(i),n}_{j+1,k}\big(u^{(i),n}_{j+1,k}\big)^2,\qquad i=m+1,\dots,p.
\end{aligned}
\end{equation*}  
One can then obtain 
\begin{equation}
\begin{aligned}
&\big(u^{(i),n+1}_{j,k}\big)^2-\big(u^{(i),n}_{j,k}\big)^2\le D^{(i),n}_{j,k}\Big(\big(u^{(i),n}_{j-1,k}\big)^2-
\big(u^{(i),n}_{j,k}\big)^2\Big),\qquad i=1,\dots,m,\\
&\big(u^{(i),n+1}_{j,k}\big)^2-\big(u^{(i),n}_{j,k}\big)^2\le D^{(i),n}_{j+1,k}\Big(\big(u^{(i),n}_{j+1,k}\big)^2-
\big(u^{(i),n}_{j,k}\big)^2\Big),\qquad i=m+1,\dots,p.
\end{aligned}
\label{c6}
\end{equation}
The estimates \eref{c6} hold provided
$$
D^{(i),n}_{j,k}+\dt\,b^{(i)}\le1,\qquad D^{(i),n}_{j+1,k}+\dt\,b^{(i)}\le1,
$$ 
which can be ensured by choosing a sufficiently small time-step $\dt$ in the upwind scheme. Continuing the proof of Theorem \ref{thm3} after
formula \eref{1.40a}, one arrives at the same estimated decay rate as in Theorem \ref{thm3}.
\begin{rmk}
It should be observed that the decay rate estimated in Theorem \ref{thm3} serves as an upper bound for the decay rate of the system
\eref{c1}--\eref{c1a} considered in this appendix; see, e.g., Example 8.
\end{rmk} 
\begin{rmk}
For simplicity, we restrict ourselves in this appendix to the case where the matrix $B$ in \eref{c1} is diagonal and positive definite. In
fact, the conclusion still holds when the matrix $B$ is positive definite, but not necessarily diagonal.
\end{rmk}

\bibliography{reference}
\bibliographystyle{siam}
\end{document}